\documentclass{article}
\usepackage{graphicx} 
\usepackage{amsmath}
\usepackage{amssymb}
\usepackage{amsthm}
\usepackage{color}

\newtheorem{theorem}{Theorem}
\newtheorem{lemma}{Lemma}
\newtheorem{algorithm}{Algorithm}
\newtheorem{conjecture}{Conjecture}

\title{Beta representations with minimal mean digit}
\author{Carl P. Dettmann, University of Bristol, UK}
\date{\today}

\begin{document}

\maketitle

\begin{abstract}
We consider the problem of minimising the mean
digit in beta representations with unrestricted digits,
equivalent to a one dimensional affine switched system
obtained from a formal limit of a stabilisable
two dimensional linear system.
We give a countable set of $\beta$ for which
the result is given by the usual (greedy)
beta expansion, an interval of values for
which it is strictly less, a numerical upper
bound, and a conditional lower bound for all
$\beta$.
\end{abstract}

\section{Introduction}
A beta representation expresses $u\in\mathbb{R}_{>0}$ to base $\beta\in\mathbb{R}_{>1}$,
\begin{equation}
    u=\sum_{j=J}^\infty\frac{d_j}{\beta^j}
\end{equation}
where $d_j\in\mathbb{Z}_{\geq 0}$ and $J\in\mathbb{Z}\cup\{-\infty\}$. 

The aim of this work is to consider beta representations with
minimal mean digit.  Typically the assumption is made that $d_j<\beta$ which implies $u\leq \beta^{1-J}\frac{\lceil\beta\rceil-1}{\beta-1}$, or a more general fixed bound~\cite{KL15} but we
will not have an a priori bound on the digits.  Of course, if $d_j>0$ for
infinitely many $j<0$ or $d_j$ grows too rapidly for $j>0$ then the series diverges and does not represent any finite $u$. Expansions using arbitrary (but finite and pre-specified) digit sets in $\mathbb{Q}(\beta)$ were considered in Ref.~\cite{KS12}.  Expansions using arbitrary non-negative digits for a finite number of non-negative powers of $\beta\in\mathbb{C}$ were considered in Ref.~\cite{Dubickas24}. 

The most widely considered beta representations, called beta expansions
and dating back to R\'enyi~\cite{Renyi57} and of active current interest~\cite{CM25,HO25,Suzuki26,Takamizo24}, are beta representations choosing the $d_j$ in a ``greedy'' manner, that is, if the remainder is
\begin{equation}\label{e:greedy}
    r_k=u-\sum_{j=J}^k\frac{d_j}{\beta^j}
\end{equation}
then $0\leq r_k<\beta^{-k}$ for all $k\geq J$.  The beta expansion exists and is unique, and satisfies $d_j<\beta$ for $j>J$.  Also $d_J=\lfloor\beta^Ju\rfloor<\beta$ iff $u<\beta^{-J}\lceil\beta\rceil$.  

There are two natural ways we can remove the bound on the digits.
Denote ${\cal D}_J\subseteq\mathbb{Z}_{\geq 0}^{\mathbb{Z}_{\geq J}}$ to be the set of digit sequences $\{d_j\}_{j=J}^\infty$ for which $d_j\in\mathbb{Z}_{\geq 0}$ and (if $J=-\infty$)
only finitely many $d_j$
are nonzero for negative $j$.  Also denote ${\cal B}_J={\cal D}_J\cap l^\infty(\mathbb{Z}_{\geq J},\mathbb{Z}_{\geq 0})$ to be the subset of ${\cal D}_J$ for which $\{d_j\}$ is bounded. In the latter case, the bound on $\{d_j\}$ is not specified a priori, but may depend on the real value $u$ represented.

The minimal mean digit is quantified as
\begin{align}\label{e:dbar}
    \bar{d}(\beta)&=\sup_{u\in\mathbb{R}_{>0}}\bar{d}_J(\beta,u)\\
    \bar{d}_J(\beta,u)&=\inf_{\{d_j\}\in{\cal D}_J}\{\limsup_{k\to\infty}\frac{1}{k}\sum_{j=J}^{k}d_j:u=\sum_{j=J}^{\infty}\frac{d_j}{\beta^j}\}\nonumber
\end{align}
where the infimum naturally excludes $\{d_j\}$ for which the series diverges.
We note that shifting $J$ is equivalent to multiplying
$u$ by powers of $\beta$, so that $\bar{d}(\beta)$ is independent
of $J$ for finite $J$.  Also, no convergent series of this form has
infinitely many $d_j>0$ for $j<0$, so that $J=-\infty$ also
gives the same result. Alternatively
\begin{align}\label{e:dbarbdd}
    \bar{d}^{(bdd)}(\beta)&=\sup_{u\in\mathbb{R}_{>0}}\bar{d}^{(bdd)}_J(\beta,u)\\
    \bar{d}^{(bdd)}_J(\beta,u)&=\inf_{\{d_j\}\in {\cal B}_J}\{\limsup_{k\to\infty}\frac{1}{k}\sum_{j=J}^{k}d_j:u=\sum_{j=J}^{\infty}\frac{d_j}{\beta^j}\}\nonumber
\end{align}
Here, since $\{d_j\}$ is bounded and $\beta>1$, the series always converges. 
Finally, we consider the beta expansion.  We denote the
mean digit in this case by 
\begin{align}\label{e:dbarbE}
    \bar{d}^{(\beta E)}(\beta)&=\sup_{u\in\mathbb{R}_{>0}}\bar{d}^{(\beta E)}_J(\beta,u)\\
    \bar{d}^{(\beta E)}_J(\beta,u)&=\limsup_{k\to\infty}\frac{1}{k}\sum_{j=J}^{k}d_{j,J}^{(\beta E)}(\beta,u)\nonumber
\end{align}
where $d_{j,J}^{(\beta E)}(\beta,u)$ is $j$th digit of the beta expansion of $u$ starting at $J$. Evidently for all $\beta$ we have
\begin{equation}
\bar{d}^{(\beta E)}(\beta)\geq\bar{d}^{(bdd)}(\beta)\geq\bar{d}(\beta)
\end{equation}
Neither $\bar{d}(\beta)$ nor $\bar{d}^{(bdd)}(\beta)$
seem to have been previously studied; the closest is Ref.~\cite{FS08}
which minimises finite sums of absolute value of (possibly negative)
bounded digits where $\beta$ is Pisot, focusing on $\rho$, $\phi$ and
$\mu_3$ (see Table~\ref{t:db1} below); see also Ref.~\cite{GS11}.
Recall that a Pisot number is an algebraic integer with all 
conjugates of modulus less than unity.  There is some literature 
directly pertinent to $\bar{d}^{(\beta E)}(\beta)$ as noted in 
Section~\ref{s:thm} below.

\begin{figure}
\centerline{\includegraphics[width=400pt]{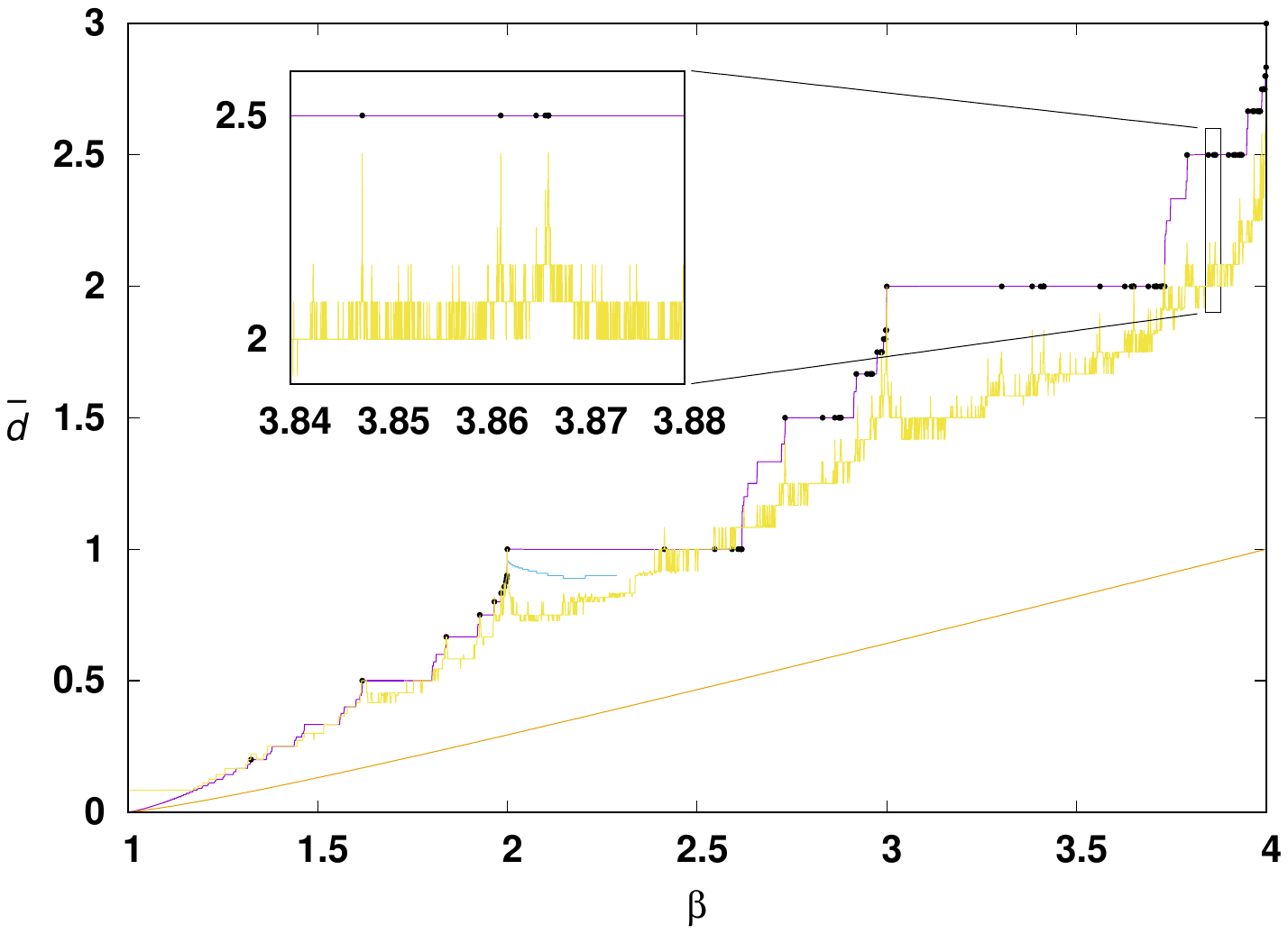}}
\caption{Upper curve: $\bar{d}^{(\beta E)}(\beta)$ (numerical).  Black dots: Elements of $\{\rho\}\cup MB$ where $\bar{d}^{(\beta E)}(\beta)=\bar{d}(\beta)$ from Theorem 1.  The small section of curve for $2<\beta<2.289$ gives the upper bound for $\bar{d}^{(bdd)}(\beta)$ from Theorem 2. The nonmonotonic yellow curve gives the numerical upper bound for $\bar{d}^{(bdd)}(\beta)$.
The lower curve is the conditional lower bound for $\bar{d}^{(bdd)}(\beta)$ from Theorem 3.
 \label{f:dbar}}
\end{figure}

This problem arose from stabilisability of switched systems.  Motivation,
relation to joint spectra of matrix sets, control theory and applications are presented in Section~\ref{s:switch} which may be
omitted at first reading, however there is a conjecture needed for the
proof of Theorem 3.  Section~\ref{s:thm} presents three theorems concerning this problem.  Theorem 1 gives a countable set of $\beta$
for which $\bar{d}^{(\beta E)}(\beta)=\bar{d}(\beta)$.  Theorem 2 gives a non-trivial upper bound on $\bar{d}^{(bdd)}(\beta)$ for $2<\beta<2.28879$ showing that $\bar{d}^{(bdd)}(\beta)<\bar{d}^{(\beta E)}(\beta)$ here. Theorem 3 gives a conditional lower bound on $\bar{d}^{(bdd)}(\beta)$ for all $\beta\in(1,\infty)$.  These theorems together with a numerical upper bound for $\bar{d}^{(bdd)}(\beta)$ are illustrated in Figure~\ref{f:dbar}.  The proofs of these theorems are given in the remaining sections.

\section*{Acknowledgements}
The author is grateful for helpful discussions with Thomas Jordan and Chenmiao Zhang, and helpful comments from an anonymous referee.  For the purpose of open access, the author has applied a Creative Commons Attribution (CC BY) licence to any Author Accepted Manuscript version arising from this submission.

\section{Stabilisability of switched systems:\\ Motivation and applications}\label{s:switch}
Joint spectral properties of sets of matrices are
notoriously challenging to characterise.  For example,
for the joint spectral radius $\rho_\infty$~\cite{RS60}
it is known that the question $\rho_\infty\leq 1$ is
undecidable~\cite{BT00}.  Here, we consider the
stabilisability radius $\tilde{\rho}$~\cite{JM17}.
Given a discrete time switched linear system
\begin{equation}\label{e:dtss}
    x(k+1)=A_{\sigma_k}x(k),\qquad x(0)=x_0
\end{equation}
defined by a set of matrices $\mathcal{M}=\{A_\sigma\}_{\sigma\in\Sigma}\subset\mathbb{R}^{n\times n}$,
$\tilde{\rho}$ characterises the ability, using a knowledge of the state of the system, to stabilise it (that is, send $x(k)$ to zero as quickly as possible) by controlling the switching signal $\{\sigma_k\}$:
\begin{align}
\tilde\rho(\mathcal{M})&=\sup_{x\in\mathbb{R}^n\setminus\{0\}}\tilde{\rho}_x(\mathcal{M})\\\nonumber
\tilde{\rho}_x(\mathcal{M})&=\inf\{\lambda\geq 0\;|\;\exists M>0,\{\sigma_k\}\; s.t. \;\|x(k)\|\leq M\lambda^k\|x\|\;\forall k\in\mathbb{Z}_{\geq 0}\}
\end{align}
where $x(k)$ is determined by~(\ref{e:dtss}) with $x_0=x$.

Many engineering control problems (mechanical, power, traffic,
etc.) involve control of multiple variables with a selection of
actions~\cite{ZZSL16}. This is also the case for
treatment regimes of
viruses and cancer, eliminating the pathogens in the presence of
mutation, drug resistance and drug toxicity~\cite{AGFH21}. Switched systems have also been used to understand
strategies of microbes that switch between dormant and active
states~\cite{BHS21}.

Consider the switched linear system~(\ref{e:dtss}) with two $2\times 2$ matrices as follows
\begin{equation}\label{e:A1A2}
\mathcal{M}(\theta,c,\beta)=\left\{A_1=\left(\begin{matrix}\cos\theta&\sin\theta\\-\sin\theta&\cos\theta\end{matrix}\right),\qquad A_2=\left(\begin{matrix}c&0\\0&\beta c\end{matrix}\right)\right\}
\end{equation}
where $\beta\geq c^{-1}>1$.  Without knowledge of the initial condition $x(0)$, it is not possible to stabilise the system, since any product of these matrices has determinant at least equal to 1.  However, given the initial condition, and for almost all values of $\theta$, the matrix $A_1$ can be used to rotate arbitrarily close to the $x_1$-axis, after which the matrix $A_2$ reduces $\|x\|$. The case $\mathcal{M}(\pi/6,1/2,4)$ was proposed in Ref.~\cite{Stanford79} and also discussed in Ref.~\cite{SU94}.  In Ref.~\cite{DJM20},
$\mathcal{M}(\theta,c,c^{-2})$ was considered, especially the special case $\mathcal{M}(\pi/4,1/2,4)$.  For the latter example, the best bounds given were
\begin{equation}
    0.707\approx \frac{1}{\sqrt{2}}\leq \tilde{\rho}(\mathcal{M}(\pi/4,1/2,4))<0.9^{1/4}\approx 0.974
\end{equation}
This is a large gap as the system is still poorly understood.

In order to better characterise these systems, we simplify
the problem further whilst retaining essential properties.  Assume that
$\theta=O(\epsilon)$ and $u=\arctan(x_2/x_1)/\theta=O(1)$ for some small $\epsilon>0$.  Then, expanding the arctan in a Maclaurin series, to linear order in $\epsilon$, the actions of the matrices $A_1$ and $A_2$ become
\begin{align}\label{e:a}
    u(k+1)&=a_{\sigma_k}(u(k))\\
    a_\sigma(u)&=\left\{\begin{array}{cc}u-1&\sigma=1\\\beta u&\sigma=2\end{array}\right.\nonumber
\end{align}
where $u\in\mathbb{R}$, $\beta\in\mathbb{R}_{>1}$ and the optimal stabilisation (choice of switching signal to make $x(k)\to 0$ as
quickly as possible) now corresponds to minimising the proportion of $a_1$ transformations.  However, the approximation is valid only if
we require $u$ to remain bounded, which is possible only if $u(0)\geq 0$ and $\sigma_k=2$ if $u(k)<1$.  

This simplification neglects two effects of the original system, firstly the nonlinear terms, of which the leading term is $O(\epsilon^2)$ with a negative coefficient in the $a_2$ equation of~(\ref{e:a}).  Secondly, whilst it can be assumed that
small rotations can eventually direct the system near the stable $x_1$-axis so that $\arctan(x_2/x_1)$ is small as desired, it is also possible that a large number of small rotations could, if necessary, rotate around to a more favourable configuration, for example, exactly onto
the $x_1$-axis.  However for very small $\epsilon$ this would be a huge penalty, so the set of initial conditions where it improves the stabilisation would be vanishingly small.  In this paper, we consider~(\ref{e:a}) as an interesting system in its own right.

Now, we show a correspondence between beta representations with
$J=0$ and the affine switched system as follows. For any bounded orbit
of the switched system~(\ref{e:a}), we note that the length of consecutive sequences of $k$ where $\sigma_k=1$ must be bounded.
Consider $\{d_j\}_{j=0}^\infty=\Phi(\{\sigma_k\}_{k=0}^\infty)$ defined as follows: Let $K_j$ for $j\in\mathbb{Z}_{\geq 1}$ be the location
of the $j$th value of $k$ for which $\sigma_k=2$. Then $d_0=K_1$ and $d_j=K_{j+1}-K_j-1$ for $j>0$ so that $\{d_j\}$ is bounded.
The inverse transformation $\Phi^{-1}$ takes bounded digit sequences $\{d_j\}$ and constructs $K_j=j-1+\sum_{i=0}^{j-1}d_i$ for $j\geq 1$.
Then $\sigma_k=2$ if $\exists j: K_j=k$ otherwise $\sigma_k=1$.

Now, consider the switched dynamics where $u(0)=\sum_{j=0}^\infty d_j\beta^{-j}$.  Whenever $\sigma_k=1$ it decreases $d_0$ by 1 and
whenever $\sigma_k=2$ (and hence $d_0=0$) it shifts the digits in the beta expansion in such a way that $\Phi(\{\sigma_j\}_{j=k}^\infty)$ always represents $u(k)$.  Thus $\Phi^{-1}$ applied to a bounded
beta representation yields a bounded solution of the switched system.
Also, for any other initial value for the switched system, the
discrepancy with the beta representation grows by a factor $\beta$
each time $\sigma_k=2$, and so is unbounded.

The above correspondence relates bounded solutions of the switched
system to beta representations with bounded digits.  As noted above,
we can require the beta representation to converge but allow the
digits not to remain bounded.  In principle we could allow the switched
system to be unbounded, but this is problematic for two reasons:
First, it could not approximate the two dimensional switched system for any finite $\theta$, and secondly, many unbounded solutions
(including negative ones) do not correspond to a beta representation.
So, we consider only the bounded switched system which has a 1:1 correspondence to beta representations with bounded digits, and we
also consider more general convergent beta representations. 

Assuming the simplifications above are valid, we can obtain the stabilisability radius $\tilde{\rho}$ of the original system~(\ref{e:A1A2}) from the beta
representation as follows.  Suppose we have a sequence of transformations
of length $T=T_1+T_2$ where $T_1$ is the number of occurrences of
$A_1$, $T_2$ is the number of occurrences of $A_2$, and the sequence ends with an $A_2$.  The matrix
$A_1$ has no effect on $\|x\|$, whilst $A_2$ multiplies it
by $c$ in the linear approximation.  The beta representation is then
\begin{equation}
    u(0)=\sum_{j=0}^{k-1}\frac{d_j}{\beta^j}+\frac{u(T)}{\beta^k}
\end{equation}
so that $T_2=k$ and
\begin{equation}
    T_1=\sum_{j=0}^{k-1}d_j=\bar{d}_kk=\bar{d}_kT_2
\end{equation}
where
\begin{equation}
    \bar{d}_k=\frac{1}{k}\sum_{j=0}^{k-1}d_j
\end{equation}
is the mean of the first $k$ digits.  Thus $T=T_1+T_2=T_2(1+\bar{d}_k)$ and
\begin{equation}  \left(\frac{\|x(T)\|}{\|x(0)\|}\right)^{1/T}=c^{T_2/T}=c^{1/(\bar{d}_k+1)}
\end{equation}
Taking the limit $T\to\infty$ we obtain
\begin{conjecture}
\begin{equation}\label{e:rhodbar}
    \lim_{\theta\to 0}\tilde{\rho}({\cal M}(\theta,c,\beta))=c^{1/(\bar{d}^{(bdd)}(\beta)+1)}
\end{equation}
\end{conjecture}
Here, the conjecture includes the statement that the limit exists, ${\cal M}(\theta,c,\beta)$ is defined in~(\ref{e:A1A2}) with conditions $\beta\geq c^{-1}>1$, and $\bar{d}^{(bdd)}(\beta)$ is given by~(\ref{e:dbarbdd}), that is, the infimum of mean digits for beta representations of the initial value $u$, using bounded digits.  The limit superior in~(\ref{e:dbarbdd}) (and similarly~(\ref{e:dbar}, \ref{e:dbarbE})) follows from the definition of $\tilde{\rho}$ leading to an upper bound on the norm of the orbit for all time. 

\begin{table}
\centerline{
    \begin{tabular}{|c|c|c|c|l|}\hline
    Symbol&Value&Pisot?&Minimal polynomial&$d_\beta(1)$\\\hline
    $\psi_6$&1.2852&N&$x^6-x^5-1$&100001\\
    $\rho$&1.3247&Y&$x^3-x-1$&10001\\
    $\psi_4$&1.3803&Y&$x^4-x^3-1$&1001\\
    $\sqrt{2}$&1.4142&N&$x^2-2$&1001000001\ldots\\
    $\psi_3$&1.4656&Y&$x^3-x^2-1$&101\\
    $\phi$&1.6180&Y&$x^2-x-1$&11\\
    $\mu_3$&1.8393&Y&$x^3-x^2-x-1$&111\\
    $\gamma_6$&2.2056&Y&$x^3-2x^2-1$&201\\
    $\gamma_5$&2.2888&N&$x^4-3x^3+x^2+x+1$&2011002001\ldots\\
    $e$&2.7183&N&&2121111212\ldots\\\hline      
    \end{tabular}}
    \caption{Beta expansion of unity for irrational numbers appearing in this paper.  We may also denote $\rho=\psi_5$ and $\phi=\psi_2=\mu_2$. \label{t:db1}}
\end{table}

\section{Results}\label{s:thm}
We now return to beta representations, giving a little more background before stating the results of this paper.
A special role is given to the beta expansion of $u=1$ for $J=1$ which is denoted $d_\beta(1)$, written $d_1d_2\ldots$ and may be finite (omitting trailing zeros) or infinite.  Note that the choice $J=1$ excludes the trivial expansion given only by $d_0=1$.  A number $\beta$ for which $d_\beta(1)$ is eventually periodic is called a Parry number, and if it is finite it is called a simple Parry number~\cite{Parry60,AMPF06}.  For example, $d_\beta(1)=\beta$ if $\beta\in\mathbb{Z}_{\geq 2}$.  Also $d_{\phi}(1)=11$ where $\phi=(1+\sqrt{5})/2$ is the golden ratio, satisfying $1=\phi^{-1}+\phi^{-2}$.  Refer to Table~1 for $d_\beta(1)$ of irrational numbers appearing in this paper.

Recall the greedy expansion, (\ref{e:greedy}). The quasi-greedy expansion of 1 instead satisfies $0<r_k\leq\beta^{-k}$ for all $k$, and is denoted $d^*_\beta(1)$.  It is always infinite, and is equal to $d_\beta(1)$ iff the latter is infinite.  For example, $d^*_\beta(1)=(\beta-1)^\infty$ if $\beta\in\mathbb{Z}_{\geq 2}$ and $d^*_{\phi}(1)=(10)^\infty$.

We define $MB$ to be the set of monotonic beta expansions, that is
\begin{equation}
MB=\{\beta>1:d_{\beta}(1)=d_1d_2\ldots,\quad d_k\geq d_{k+1} \quad\forall k\in\mathbb{Z}_{\geq 1}\}
\end{equation}
This incorporates both classes considered by Frougny and Solomyak~\cite{FS92} in their Theorem 2 (finite case) and Theorem
3 (infinite case).  In both cases it was shown that elements of
$MB$ are Pisot numbers.  We note $\mathbb{Z}_{\geq 2}\subset MB$.

The most obvious approach to minimise the mean digit is the
greedy strategy, choosing digits as large as possible or equivalently choosing the transformation $a_1$ when possible.  This gives the beta expansion of the initial point $u$ with $J=0$.

We now recall the definition of $\bar{d}^{(\beta E)}$,~(\ref{e:dbarbE}).  Digit frequencies of beta expansions have been studied in Ref.~\cite{BCH16}.  In particular, for Lebesgue almost every $\beta$,
digit frequencies form a polytope with rational vertices, for which $\bar{d}^{(\beta E)}$ gives the largest mean digit.  A given polytope exists for a finite interval in $\beta$,
thus $\bar{d}^{(\beta E)}$ is almost everywhere locally constant as a function of $\beta$, as
illustrated in Figure~\ref{f:dbar}.  The Hausdorff dimension of sets where the mean digit does not exist ($\liminf$ strictly less than $\limsup$) for beta
expansions was considered in Ref.~\cite{LL16}.

As an aside, Ref.~\cite{DDKL12} has an alternative notion of ``each
step'' optimal beta representations for which only a beta expansion
can be optimal, and is so for all $u$ only for $d_\beta(1)=m^np$
where $1\leq p\leq m$, $0\leq n<\infty$ and $\beta>1$, that is,
a subset of MB.

Beta expansions have the property that no sequence of digits
(after $d_J$) can be greater lexicographically than $d_\beta(1)$~\cite{Parry60}.
Thus it is (for almost all $\beta$) a finite (and short) process
to calculate $\bar{d}^{(\beta E)}$.
For example, for $\beta=e\approx 2.7183$, $d_j<e$ so the maximum digit is 2.
Looking at $d_e(1)$, the sequence $22$ is disallowed, and the sequence
$212$ is allowed but only followed by $111$, $102$, $021$ or other sequences with lower mean.  A repeating sequence of $211$ is allowed.  So, $\bar{d}^{(\beta E)}=4/3$.

The numerical computation used to generate the $\bar{d}^{(\beta E)}(\beta)$ curve in Figure~\ref{f:dbar} used digit sequences up to length
20 for $2<\beta\leq 4$ where the output was always found to be rational with denominator at most 7, digit sequences up to length 30 for
$1.088\leq\beta\leq 2$, and for $1<\beta<1.088$ the asymptotic formula
\begin{equation}
    \bar{d}^{(\beta E)}(\beta)\sim \frac{\ln\beta}{\ln(\beta/(\beta-1))},\qquad \beta\to1^+
\end{equation}
which follows by considering the sequence of $\beta$ such that $d_\beta(1)=10^{n-2}1$ or equivalently $d^*_\beta(1)=(10^{n-1})^\infty$
and which have $\bar{d}^{(\beta E)}(\beta)=n^{-1}$, and noting that
$\bar{d}^{(\beta E)}(\beta)$ is monotonic.  The locally constant
property can make the computation more efficient:
If $\bar{d}^{(\beta E)}(\beta_1)=\bar{d}^{(\beta E)}(\beta_2)$
then $\bar{d}^{(\beta E)}(\beta)=\bar{d}^{(\beta E)}(\beta_1)$ for all
$\beta\in[\beta_1,\beta_2]$.

If $\beta$ is an integer, it is clear that if any digits are greater than $\beta-1$, multiples of $\beta$ may be removed with equivalent value added at lower $j$ without increasing the digit mean.
Also, some real numbers have a full density of $\beta-1$ digits. So the beta representation is optimal and $\bar{d}(\beta)=\beta-1$ in this case. But for all but a countable set of $\beta$ the greedy strategy is not always optimal:

\begin{theorem}\label{th:bO}
    Let $\beta\in\mathbb{R}_{>1}$, $J\in\mathbb{Z}\cup\{-\infty\}$.\\
    (a) If $\beta\in\{\rho\}\cup MB$, then $\bar{d}(\beta)=\bar{d}^{(\beta E)}(\beta)$.\\
    (b) If $\beta\not\in\{\rho\}\cup MB$ then there exists a positive integer $m$ such that the beta representation has a larger digit sum than $m$ and furthermore, if\\
    (i) the beta representation of $m$ is finite, or\\
    (ii) there is an integer $m'$ so that $\bar{d}_J^{(\beta E)}(\beta,m')>0$, or\\
    (iii) $d_\beta(1)$ is finite or has an upper bound on the number of consecutive zeros, then $\exists u\in\mathbb{R}_{>0}$ such that $\bar{d}_J(\beta,u)<\bar{d}_J^{(\beta E)}(\beta,u)$.
\end{theorem}
Here, $\rho\approx 1.3247$ is the real rot of $x^3-x-1$ (the minimal Pisot number); see Table~\ref{t:db1}.  The set $MB$ discussed above includes the Fibonacci (golden), tribonacci and higher multinacci numbers $\mu_k$ (with $k\geq 2$) 
that satisfy $x^k=\frac{x^k-1}{x-1}$ with $x>1$.
We have $2-\mu_k=\mu_k^{-k}\sim 2^{-k}$ as $k\to\infty$, and $\bar{d}^{(\beta E)}(\mu_k)=1-k^{-1}$ which is thus not H\"older
continuous at $\beta=2$.
Theorem 1(a) then implies that the functions $\bar{d}^{(bdd)}(\beta)$ and $\bar{d}(\beta)$ are also not H\"older continuous at $\beta=2$.

In Theorem 1(b), if $\beta>\sqrt{2}$ then we may take $m=\lfloor\beta\rfloor+1$.  Beta expansions of integers to Pisot
bases are finite or eventually periodic~\cite{Schmidt80}, so all
Pisot numbers are included in (i) or (ii).  Condition (ii) is valid
apart from an extremely small set of $\beta$, presumably of zero Hausdorff dimension (see~\cite{LL16} for a related problem).  Condition (iii) is presumably a set of full Hausdorff dimension but zero Lebesgue measure (similar to continued fractions with bounded digits), but may be useful if information is known about $d_\beta(1)$.
For example $\psi_6\approx1.2852$ has finite $d_\beta(1)$, is not
Pisot, gives an apparently aperiodic beta expansion
$2=100.0^410^610^510^510^610^81\ldots$ and so is not clearly
covered by conditions (i) or (ii).  Here, $\psi_k$ (for $k\geq 2$) 
denotes the positive root of $x^k-x^{k-1}-1$, generalising the 
supergolden ratio $\psi_3$; see Table 1.

Theorem 1(b) says that for almost all $\beta$ there are some $u$ for which
the greedy algorithm is not optimal, however this is not sufficient to show
$\bar{d}(\beta)<\bar{d}^{(\beta E)}(\beta)$. 
Now, we give an explicit non-trivial upper bound for $\bar{d}^{(bdd)}(\beta)$ and hence for $\bar{d}(\beta)$, for all $\beta$ in an interval:
\begin{theorem}
    The equation $4x^{k-2}=\frac{x^{k}-1}{x-1}$, that is, the value of the base $\beta$ in which $040^{k-2}=1^k$, has
    a unique real solution, denoted $\gamma_k$, for $x>2$ and each $k\geq 5$. Then for $2<\beta<\gamma_5\approx 2.28879$, $\bar{d}^{(\beta E)}(\beta)=1$ but\\
    (a) For $\gamma_{6}\leq \beta\leq \gamma_{5}$,
        $\bar{d}^{(bdd)}(\beta)\leq \frac{9}{10}$.\\
    (b) For $\gamma_{k+1}\leq\beta\leq\gamma_{k}$
        with $k\geq 6$, $\bar{d}^{(bdd)}(\beta)\leq \frac{k+2}{k+3}$.
\end{theorem}
Given that $\bar{d}(2)=1$ from Theorem 1(a), this shows that $\bar{d}^{(bdd)}(\beta)$ and $\bar{d}(\beta)$ are not increasing functions.
Moreover, noting from~(\ref{e:proof2})
below that $\gamma_k-2$ is exponentially
decreasing as $k\to\infty$, we find also from $\beta>2$ that $\bar{d}^{(bdd)}(\beta)$ and $\bar{d}(\beta)$ are not H\"older continuous at $\beta=2$.
For $\beta=\gamma_k$ we can choose
the best bound, so $\bar{d}(\gamma_5)\leq \frac{9}{10}$,
$\bar{d}(\gamma_6)\leq \frac{8}{9}$, $\bar{d}(\gamma_k)\leq \frac{k+1}{k+2}$ for $k\geq 7$.

In Figure~\ref{f:dbar} we also provide a numerical
upper bound for $\bar{d}^{(bdd)}(\beta)$.
This is obtained as follows.  For
each $\beta$, consider digit sequences of length
$k$ for $k\leq 12$.  These are enumerated in
increasing order of digit sum.  When the beta
representations corresponding to these digit
sequences cover the unit interval with
gap at most $\beta^{-k}$ then we can construct
a beta representation for any $u\in[0,1)$ using
these sequences, $k$ digits at a time.
With $d_0=\lfloor u\rfloor$,
this then gives the same mean digit for
$u\in\mathbb{R}_{>0}$ and we have
an upper bound for $\bar{d}^{(bdd)}(\beta)$.

An efficient method of checking coverage is to
divide the unit interval into bins of size
$\beta^{-k}$ and record the minimum and maximum
value in each bin.  If all bins are occupied and
the gap between the maximum of one bin and the
minimum of the next is at most $\beta^{-k}$
then coverage is obtained.  If the digit
sequences vary by one in the final digit, the
gap is exactly $\beta^{-k}$.  Otherwise a small
tolerance on the gap size is needed to avoid round off errors.
This is taken to be $10^{-14}$, somewhat greater than the
limit of double precision ($2^{-53}\approx 10^{-16}$).
In this calculation the propagation of errors was not
checked rigorously.

The number of digit sequences $\{d_j\}_{1\leq j\leq k}$ with $d_j\geq 0$ and $\sum_jd_j\leq\lfloor k\bar{d}\rfloor$ is given by the stars and bars formula
\begin{equation}\label{e:k=12}
\left(\begin{array}{c}\lfloor k\bar{d}\rfloor+k\\k\end{array}\right)=\left(\begin{array}{c}48\\12\end{array}\right)\approx 7\times 10^{10}
\end{equation}
for $k=12$ and $\bar{d}=3$.

It is seen in Figure~\ref{f:dbar} that the numerical upper
bound for 
$\bar{d}^{(bdd)}(\beta)$ is mostly below $\bar{d}^{(\beta E)}$,
however it rises near the points in $\rho\cup MB$
as required by Theorem 1.  The inset to the figure
is given to illustrate this on a smaller scale.
It is found that increasing $k$ usually leads to
an improvement of the bound, so that for the vast
majority of values of $\beta$ in the figure, the best bound
was for $k=12$ and for almost all the rest, $k=11$.  This
also suggests that modest improvement is possible by
intensive calculations exceeding $k=12$, noting~(\ref{e:k=12}). 

Finally, we give a conditional lower bound for $\bar{d}^{(bdd)}(\beta)$.  Let us define $F(x)=\frac{(x+1)^{x+1}}{x^x}$ for $x\in\mathbb{R}_{\geq 0
}$, where as usual $0^0=1$.  This is a strictly increasing function mapping $[0,\infty)$ to $[1,\infty)$ with inverse $F^{-1}(x)$. We have
\begin{theorem} 
Assuming Conjecture 1,
\begin{align}
\bar{d}^{(bdd)}(\beta)&\geq F^{-1}(\beta)
\end{align}
\end{theorem}
In addition to the proof given
in Section~\ref{s:thm3}, there is a naive argument for this
result, namely, asserting that we require $\beta^k$ sequences
of length $k$, comparing this with~(\ref{e:k=12}) and
taking $k\to\infty$.  However, the mean digit constraint
applies only at large $k$, and not for each fixed $k$.

It would be good to prove Theorem 3 without relying on Conjecture 1, and to resolve the conjecture. From a practical point of view, the
conjecture, together with Theorem 1, suggest that systems of the
form (\ref{e:A1A2}) are more difficult to stabilise if $\beta\in\{\rho\}\cup MB$.  In addition, it would be interesting to
find $\beta$ such that $\bar{d}(\beta)<\bar{d}^{(bdd)}(\beta)$ or show
this does not occur. Looking at Figure~\ref{f:dbar},  we see that as for the $2\times 2$ switched system, there is still a large gap between upper and lower bounds, hence much scope for further work.

\section{Proof of Theorem 1} 

Now we propose an algorithm to reduce an arbitrary beta representation
to the unrestricted ($J=-\infty$) beta expansion.

\begin{algorithm}
Input: Value $\beta\in\mathbb{R}_{>1}$.  If the beta  expansion of unity is finite, $d_{\beta}(1)=d_1d_2\ldots d_k$, the set of disallowed words is
$\{(d_1+1),d_1(d_2+1),\ldots,d_1d_2\ldots d_{k-2}(d_{k-1}+1),d_1\ldots d_{k-1}d_k\}$.  If it is infinite, $d_{\beta}(1)=d_1d_2\ldots$, the set of disallowed words is $\{(d_1+1),d_1(d_2+1),\ldots,d_1\ldots d_{k-1}(d_k+1),\ldots,d_\beta(1)\}$.\\
Input: Arbitrary beta representation
\[ u=\sum_{j=-\infty}^\infty\frac{d_j}{\beta^j}\]
Step 1: Calculate the beta expansions of the disallowed words.\\
Step 2: Find the lowest value of $j$ at which the representation is disallowed, or for which the digit that is too high exceeds that of the disallowed word.
If no such value exists, stop.  If more than one disallowed word corresponds to the lowest $j$, choose the first on the above lists.\\
Step 3: Subtract the word found in step 2, and add its beta expansion.
Output the current beta representation and go to step 2.
\end{algorithm}
For example, if $\beta=\phi$, the golden ratio, $d_{\phi}(1)=11$ and the minimal disallowed words are $2=10.01$ and $11=100$.  Applying
the algorithm, we find $5=13.01=21.02=101.12=110.02=1000.02=1000.1001$.  Note
that in the second equality, both disallowed words occur at $j=0$ and we have chosen the first.  This rule is only for definiteness; the algorithm works for an arbitrary choice.

This algorithm may require manipulation of infinite sequences of disallowed
words and an infinite beta representation.  Thus a practical implementation
requires suitable truncation of these sequences.  We are not concerned with
a practical implementation here, though, only the output, which is a
sequence of beta representations of $u$ that may terminate.  What can we
say about the outcome of Algorithm 1?  We have

\begin{lemma}
Where the beta expansion of each disallowed word has less or equal digit sum than the word itself, the output of Algorithm 1 stops at, or converges to, the beta expansion of $u$, which has less or equal digit mean than the initial beta representation.
\end{lemma}

\proof
Given a fixed integer $k$, there are a finite number of combinations of digits $\{d_j\}_{j<k}$ in beta representations of $u$, since  $\sum_jd_j\beta^{j}\leq u$ and hence $d_j\leq u\beta^j$, including $d_j=0$ for $j<-\ln u/\ln \beta$.  Each replacement of a disallowed word gives a greater lexicographic outcome, so the representation for
$j<k$ converges in a finite number of steps.  Thus, the algorithm either stops at, or converges to, the beta expansion of $u$.

We write the initial beta representation as
$u=\sum_{j=J}^\infty d_j\beta^{-j}$
with $\bar{d}=\limsup_{k\to\infty}k^{-1}\sum_{j=J}^kd_j$.
Thus for all $\varepsilon>0$ there exists $K$ such that 
\begin{equation}\label{e:dbound}
\sum_{j=J}^kd_j<k(\bar{d}+\varepsilon)
\end{equation}
for all $k>K$. If we take the finite representation
$d_Jd_{J+1}\ldots d_k$
for any $k>K$ and apply the algorithm, at each step the digit
sum is nonincreasing. Hence the algorithm stops after a finite
number of steps, yielding the beta representation 
$d'_Jd'_{J+1}\ldots$, with digit sum less than or equal to 
$\sum_{j=J}^kd_j$.  In particular 
$\sum_{j=J}^kd'_j\leq\sum_{j=J}^kd_j$.  

We can thus write
\begin{equation}
    u=\sum_{j=J}^k\frac{d'_j}{\beta^j}+\sum_{j=k+1}^\infty
    \frac{d_j+d'_j}{\beta^j}
\end{equation}
Given (\ref{e:dbound}) and that $d'_j$ are from a beta
expansion, the second term satisfies
\begin{equation}
  \sum_{j=k+1}^\infty
    \frac{d'_j+d_j}{\beta^j} \leq
    \frac{(k+1)(\bar{d}+\varepsilon)}{\beta^{k+1}}+
    \sum_{j=k+2}^\infty\frac{\bar{d}+\varepsilon}{\beta^j}
    +\sum_{j=k+1}^\infty\frac{\lfloor\beta\rfloor}{\beta^j}
\end{equation}
where the first two terms give the largest possible value
from the $d_j$ contribution and the third term an upper bound
on the $d'_j$ contribution.  The beta expansion of
this term thus affects $O(\log k)$ digits for $j\leq k$.
Thus the beta expansion of $u$ has digit sum for
$j\leq k$ is at most $k(\bar{d}+\varepsilon)+O(\log k)$.
The algorithm applied to the representation of $u$
converges to this beta expansion, and we can take $k\to\infty$ and $\varepsilon\to 0$, so the value of $\bar{d}$ does not increase
as claimed. \qed

The above example, finding the beta representation of 5 base $\phi$, illustrates the lemma.

\vspace{10pt}

{\bf Proof of Theorem~\ref{th:bO}:}\\
(a) We show that all beta representations for $\beta\in\{\rho\}\cup MB$ can be reduced to beta expansions without increasing their digit sum, using Algorithm 1 and Lemma 1. Thus, $\bar{d}(\beta)=\bar{d}^{(\beta E)}(\beta)$.

For $\beta=\rho$, we have $d_{\beta}(1)=10001$ so the disallowed sequences are: $2=100.00001$, $11=1000$, $101=1000.001$, $1001=10000.00001$, $10001=100000$.

For $\beta\in MB$, we have $d_\beta(1)=d_1d_2\ldots$ with $d_k\geq d_{k+1}$ for $k\geq 1$.  Assume for now that $d_\beta(1)$ is infinite. The disallowed sequences are $d_1d_2\ldots(d_k+1)$ for
$k\geq 1$.  We have
\begin{align}\nonumber
    0.d_1d_2\ldots(d_k+1)&=0.d_1d_2\ldots d_kd_1d_2\ldots\\
    &=1.0^k(d_1-d_{k+1})(d_2-d_{k+2})\ldots\label{e:re}
\end{align}
where $0^k$ is a sequence of $k$ zeros.  Since the sequence $d_k$
is monotonic, all the digits are non-negative.  Furthermore, they form a sequence lower lexicographically than $d_\beta(1)$ and so are the beta expansion.  When $d_\beta(1)$ is finite, the above argument holds, and for the final disallowed sequence we have
\begin{align}\nonumber
    0.d_1d_2\ldots d_k=1
\end{align}
For either $\beta=\rho$ or $\beta\in MB$, the sum of the digits has not been increased by
reducing to the beta expansion, as required.\\
(b) The proof is in two parts.  First, we show that for any $\beta\not\in MB$ there is an integer for which the beta expansion has greater or infinite digit sum.  Then, for the cases given in the theorem, we use this to find $u$ for which the mean digit is strictly larger.

If $\beta>\sqrt{2}$, write $d_1=\lfloor\beta\rfloor$ and suppose that a digit $d_1+1$ can be reduced to a beta expansion without increasing the digit sum.  This then takes the form
\begin{equation}\label{e:d+1}
    (d_1+1)=10.b_1b_2\ldots
\end{equation}
for non-negative integers $b_k$.  In order not to increase the digit sum, we have
\begin{equation}\label{e:sumb}
    d_1\geq\sum_{k=1}^\infty b_k
\end{equation}
and in particular that all but a finite number of $b_k$ are zero.

Equation~(\ref{e:re}) above is valid for arbitrary $\beta$, though some of the digits might be negative. Putting $k=1$ and comparing with~(\ref{e:d+1}) allows us to
express the $d_j$ in terms of the $b_j$ as
\begin{equation}
    d_{\beta}(1)=d_1(d_1-b_1)(d_1-b_1-b_2)\ldots
\end{equation}
which has only non-negative digits due to~(\ref{e:sumb}), and
is non-increasing.  Thus $\beta\in MB$.

Now, consider the case $\beta\leq \sqrt{2}$.  Thus to base $\beta$, $2\geq 100$ and the above argument does not apply.
However, if all integers have finite beta expansions, $\beta$ must be Pisot~\cite{FS92}.  There are only two Pisot numbers in
$(1,\sqrt{2}]$, $\rho$ (considered in part (a)), and $\psi_4\approx 1.3803$, see Table 1 and Ref.~\cite{DP55}.  However to base $\psi_4$
\begin{equation}\label{e:chi2}
   2=100.0^3(0^41)^\infty
\end{equation}
that is, the beta expansion is not finite and has infinite digit sum.

Having shown there is a positive integer $m$ for which the digit sum of the beta representation is greater if $\beta\not\in MB$, we now consider the mean digit of beta representations.  If as in~(\ref{e:chi2}) the mean digit is
positive, we may take $u=m$.  This covers condition (ii).  So, now assume that $\bar{d}^{(\beta E)}_{J}(\beta,m)=0$, either because the beta expansion is finite, or because the sequence of non-zero digits has zero density.  The idea of the proof is that such a
sequence has long stretches of zeros that we can fill by adding terms
of the form $m\beta^{-kM}$ for fixed $M$ and different $k$, so that
$u=m\sum_{k\in{\cal K}}\beta^{-kM}$ for a suitable ${\cal K}\subseteq\mathbb{Z}_{\geq 0}$.

If the beta representation of $m$ is finite (condition (i)), then we may take ${\cal K}=\mathbb{Z}_{\geq 0}$ leading to a periodic beta representation. Here, we need to choose $M$ large enough to avoid a carry.

We are left with condition (iii).  In this case, we set $\Delta-1$ to be the upper bound on the number of consecutive zeros in $d_\beta(1)$,
excluding trailing zeros if it is finite.
The sum of the digits in the beta expansion is greater than $m$, so let
$\tilde{M}\in\mathbb{Z}_{>0}$ be a length at which the sum is at least $m+1$, more precisely we have in the beta representation of $m$,
\begin{equation}
    \sum_{j=J}^{J+\tilde{M}-1}d_j\geq m+1
\end{equation}
Then, let $\tilde{\delta}>0$ be such that the beta representation of
$m+\delta$ is the same for $j\leq J+\tilde{M}-1$ for all $0<\delta<\tilde{\delta}$.  
We can set $M=\tilde{M}+\Delta+\lceil \frac{-\ln \tilde{\delta}}{\ln\beta}\rceil$.  Now consider blocks of size $M$,
of which the $k$th block is $J-\Delta+kM\leq j\leq J-\Delta+(k+1)M-1$.
A block is called empty if $d_j=0$ for all $j$ in the block.
We execute the following algorithm

\begin{algorithm}
Input: Values of $\beta$, $m$, $M$.\\
Step 1: Set ${\cal K}=\{0\}$\\
Step 2: Calculate the beta representation of
$u=m\sum_{k\in{\cal K}}\beta^{-kM}$\\
Step 3: If the representation has positive mean digit, that is, if
\[\limsup_{k\to\infty}\frac{1}{k}\sum_{j=J}^kd_j>0\]
then stop.\\
Step 4: Locate the empty block at lowest $k\geq 0$, append this
value of $k$ to $\cal K$, and go to step 2.
\end{algorithm}

When an integer $k$ is appended to $\cal K$ (``the operation''),
the beta representation of $m$
is added to the empty $k$th block.  The buffer of $\Delta$ zeros for $J-\Delta+kM\leq j<J+kM$ prevents any effect to lower $k$ blocks. Prior to the operation, the $k$th block is empty for $j\leq J-\Delta+(k+1)M-1=J+kM+\tilde{M}+\lceil \frac{-\ln \tilde{\delta}}{\ln\beta}\rceil$.
Thus, any nonzero digits for larger $j$ contribute less than $\beta^{-J-kM-\lceil\frac{-\ln \tilde{\delta}}{\ln\beta}\rceil}\leq\tilde{\delta}\beta^{-J-kM}$ and hence do not affect the $\tilde{M}$ digits starting from $J+kM$.
Thus after the operation, the $k$th block digit sum is at least $m+1$.

If Algorithm 2 stops, we have immediately that $\bar{d}(\beta,u)=0$ and
$\bar{d}^{(\beta E)}(\beta,u)>0$.  Otherwise, we take the limiting value
of $u$, for which the beta expansion has no empty blocks.

Now for each of the conditions (i), (ii) and (iii) we bound the mean digits.
Define $\bar{d}_{\cal K}=\limsup_{k\to\infty}\frac{1}{k}\#\{{\cal K}\cap[0,k-1]\}$ to be the upper density of $\cal K$.  The $k$th block of the beta expansion has digit sum at at least 1 if $k\not\in{\cal K}$ and
at least $m+1$ if $k\in{\cal K}$, whilst the representation containing
$m$ at $j=J+kM$ for $k\in{\cal K}$ has digit sums 0 and $m$ respectively.
Thus
\begin{equation}
    \bar{d}_J(\beta,u)\leq \frac{m\bar{d}_{\cal K}}{M}<\frac{1+m\bar{d}_{\cal K}}{M}\leq \bar{d}_J^{(\beta E)}(\beta,u)
\end{equation}
\qed

\section{Proof of Theorem 2}
Multiplying the equation for $\gamma_k$ by $x-1$
we find
\begin{equation}\label{e:proof2}
    G(k,x)\equiv x^{k-2}(x-2)^2=1
\end{equation}
Temporarily considering $k$ as a real variable, $G(k,x)$ for $x>2$ and
$k\geq 5$ is continuous and increasing in both $k$ and $x$,  $\lim_{x\to\infty}G(k,x)=\infty$, $\lim_{k\to\infty}G(k,x)=\infty$ and
$\lim_{x\to 2}G(k,x)=0$. Thus for fixed $k$ there is exactly one
solution $\gamma_k$ which decreases with $k$, and
$\lim_{k\to\infty}\gamma_k=2$.

For each $k\geq 5$ and $\gamma_{k+1}\leq\beta\leq\gamma_k$, we have in base $\beta$
the inequalities $040^{k-1}\leq 1^{k+1}\leq 040^{k-2}1$.
Thus a sequence of $k+1$ consecutive 1s followed by any
beta expansion may be replaced by $040^{k-1}$ or $040^{k-2}1$
followed by some beta expansion.

Now, for $\beta=\gamma_5$ we have $d_\beta(1)=20110\ldots$,
while for $\beta=\gamma_6$ we have $d_\beta(1)=201$.  For
$2\leq \beta\leq \gamma_5$, $\bar{d}^{(\beta E)}(\beta)=1$
since the beta expansion can have the sequence $1^\infty$
but each $2$ is followed by a $0$.

For any $u\in[0,1)$, construct its beta expansion.  Where there are at least $k+1$ consecutive 1s, use the above replacement, which gives a sequence of $k+1$ digits with mean at most $\frac{5}{k+1}$.  Alternatively, there could be $\leq k$ consecutive 1s followed by 0, which is
a sequence with mean at most $\frac{k}{k+1}$.  Finally, there could be $\leq k$ consecutive 1s followed by a 2.
In this case, continue the sequence until the mean drops below 1.  For $k=5$ the worst case is $1^k20110$ with digit mean $\frac{9}{10}$, whilst for $k\geq 6$ the worst case is $1^k200$ with digit mean $\frac{k+2}{k+3}$.  The remainder of the expansion is also a beta expansion, so this process can be used iteratively, giving the bounds on the digit means for the entire beta representation. \qed

\section{Proof of Theorem 3}\label{s:thm3}
We use the lower bound given in Theorem 3 of Ref.~\cite{DJM20}, setting $c=\beta^{-1}$. In the notation
of that paper, $m=n=2$, $\delta_1=1$, $\delta_2=\beta^{-1}$, $\Delta_1=1$, $\Delta_2=\beta^{-1}$, $\nu=(\nu_1,1-\nu_1)$
with $\nu_1\in[0,1]$,
$\bar{\nu}=(\frac{1}{1+\beta^{-1}},\frac{\beta^{-1}}{1+\beta^{-1}})$,
\begin{align}
    \Psi(\nu)&=\nu_1\log\nu_1+(1-\nu_1)\log(\beta(1-\nu_1))\\
    \Psi(\bar{\nu})&=-\log(1+\beta^{-1})<0
\end{align}
Thus, part (b) of the theorem applies and a solution of
$\Psi(\nu)=0$ exists.  We write it in the form
\begin{equation}\label{e:beta}
    \beta=\frac{\left(\frac{\nu_1}{1-\nu_1}+1\right)^{\frac{\nu_1}{1-\nu_1}+1}}{\left(\frac{\nu_1}{1-\nu_1}\right)^{\frac{\nu_1}{1-\nu_1}}}
\end{equation}
and note that due to monotonicity, it is a unique solution
in $\nu_1$.  Then the theorem states that
\begin{equation}
    \tilde{\rho}\geq \delta_1^{\nu_1}\delta_2^{\nu_2}=\beta^{-(1-\nu_1)}
\end{equation}
Applying the conjecture and substituting $c=\beta^{-1}$, we have
\begin{equation}
    \bar{d}^{(bdd)}(\beta)\geq\frac{\nu_1}{1-\nu_1}
\end{equation}
which leads to the result by comparison with~(\ref{e:beta}).
\qed

\bibliographystyle{custom3}
\bibliography{switch}

\end{document}